\documentclass[a4paper,12pt]{article}

\newtheorem{theorem}{Theorem}

\usepackage[latin9]{inputenc}
\usepackage{units}
\usepackage{amsmath}
\usepackage{amssymb}
\usepackage{graphicx}
\usepackage{color}
\usepackage{tabls}

\def\bR{\mathbb{R}}
\def\ba{\boldsymbol{\alpha}}
\def\bRn#1{\mathbb{R}^{#1}}

\def\Rxn#1{\mathbb{R}_{#1}[x]}
\def\Rxxn#1{\mathbb{R}_{#1}[\mathbf{x}]}

\def\Pxn#1{\boldsymbol{\psi}_{#1}(x)}

\def\Pxnt#1{\boldsymbol{\psi}_{#1}^{\top}(x)}
\def\Pxxn#1{\boldsymbol{\psi}_{#1}(\mathbf{x})}
\def\Pxxnt#1{\boldsymbol{\psi}_{#1}^{\top}(\mathbf{x})}

\newtheorem{problem}{Problem}
\makeatletter
\providecommand{\tabularnewline}{\\}
\@ifundefined{showcaptionsetup}{}{%
 \PassOptionsToPackage{caption=false}{subfig}}
\usepackage{subfig}
\makeatother

\author{Jan Heller${}^1$, Didier Henrion${}^{2,3,1}$, Tom\'{a}\v{s} Pajdla${}^1$}

\title{Stable radial distortion calibration by polynomial matrix inequalities programming} 

\begin{document}

\date{}
\maketitle

\footnotetext[1]{Faculty of Electrical Engineering, Czech Technical University in Prague, CZ-16627 Praha~6, Technick\'{a}~2, Czech Republic.}  
\footnotetext[2]{CNRS, LAAS, 7~avenue du colonel Roche, F-31400 Toulouse, France.}
\footnotetext[3]{Universit\'e de Toulouse; F-31400 Toulouse; France.}

\begin{abstract}               
Polynomial and rational functions are the number one choice when it
comes to modeling of radial distortion of lenses. However, several
extrapolation and numerical issues may arise while using these functions
that have not been covered by the literature much so far. In this
paper, we identify these problems and show how to deal with them by
enforcing nonnegativity of certain polynomials. Further, we show how
to model these nonnegativities using polynomial matrix inequalities
(PMI) and how to estimate the radial distortion parameters subject
to PMI constraints using semidefinite programming (SDP). Finally,
we suggest several approaches on how to incorporate the proposed method
into the overall camera calibration procedure.
\end{abstract}

\section{Introduction}

Radial distortion modeling is the most important non-linear part of
the camera calibration process~\cite{Hartley2003b}. The first works
on the topic came from the photogrammetric community~\cite{Brown66,Brown71,Slama80}.
Since then, a plethora of models has been suggested in the literature~\cite{Sturm11}.
Among the proposed models, the ones based on polynomial and rational
functions are the most popular. This popularity undoubtedly stems
from the fact that these function are easily manipulated and yet provide
sufficient fitting power for wide range or distortions. Unfortunately,
the extrapolation qualities of polynomials can be quite unpredictable
in situations where little or no data is available. However, even
if data points are missing, the overall shape of the distortion is
known \emph{a priori} in many calibration scenarios, \emph{e.g.},
the lens introduces barrel or pincushion distortions. Based on such
\emph{a priori} information, the shape of the polynomial and rational
distortion functions can be controlled by enforcing nonnegativity
of certain polynomials. For example, in the case of pincushion distortion
we can accomplish the desired shape by enforcing nonnegativity of
the first and the second derivatives of the distortion function on
the whole field of view of the camera.

In this paper, we propose a radial distortion calibration procedure
where a polynomial cost function, \emph{e.g.}, reprojection error,
is minimized subject to such shape constraints. This shape optimization
procedure is designed to stabilize the shape of the distortion function.
It is based on polynomial matrix inequalities (PMI) programming and
can be easily incorporated into an existing camera calibration procedure.

In Section~\ref{sec:camera_radial_distortion}, we formally introduce
the radial distortion function and present several extrapolation issues
arising while using polynomial and rational distortion models. Next,
in Section~\ref{sec:theory} we provide a minimal theoretical background
needed for our shape stabilization approach. In Section~\ref{sec:shape-optimization},
we demonstrate the proposed method on three types of radial distortion
shapes and models and show how to incorporate the method into an overall
camera calibration procedure. Finally, in Section~\ref{sec:experiments}
we experimentally validate our approach and show that the method guarantees
the correct shape of a distortion function without compromising the
quality of the overall camera calibration as measured by the reprojection
error.

\section{Camera Radial Distortion\label{sec:camera_radial_distortion}}

Let us suppose that a set of scene points $\mathbf{X}_{i}\in\bRn3$,
$i=1,\dots,n$ is observed by a camera. If $\mathtt{R}\in SO(3)$,
$\mathbf{t}\in\bRn3$ are the camera extrinsic parameters, a scene
point $\mathbf{X}_{i}$ gets projected into an image point $(x_{i},y_{i},1)^{\top}$:
\[
\lambda_{i}(x_{i},y_{i},1)^{\top}=\mathtt{R}\mathbf{X}_{i}+\mathbf{t},\,\lambda_{i}\in\bR.
\]
In reality, some amount of radial distortion is always present and
the camera observes a point $(\hat{x}_{i},\hat{y}_{i},1)^{\top}$
which does not coincide with the ideal (and unobservable) point $(x_{i},y_{i},1)^{\top}$.
In pixel coordinates, the camera observes a point $\mathtt{K}(\hat{x}_{i},\hat{y}_{i},1)^{\top}$,
where $\mathtt{K}\in\bRn{3\times3}$ is the matrix of intrinsic camera
parameters, the so-called calibration matrix.\emph{ Radial distortion
function} $L\colon\bR\rightarrow\bR$ is a function of radius $r=\sqrt{x_{i}^{2}+y_{i}^{2}}$
that models the radial displacement of the ideal image point position
from the center of the radial distortion as
\begin{equation}
\left(\begin{array}{c}
\hat{x}_{i}\\
\hat{y}_{i}
\end{array}\right)=L(r)\left(\begin{array}{c}
x_{i}\\
y_{i}
\end{array}\right).\label{eq:rdfunc}
\end{equation}
The function $L(r)$ is only defined for $r>0$ and $L(0)=1$, $L(r)>0$.
For the purposes of demonstration of the proposed shape optimization
procedure, we will use $L(r)$ defined as follows
\begin{equation}
L(r)=\frac{f(r)}{g(r)}=\frac{1+k_{1}r+k_{2}r^{2}+k_{3}r^{3}}{1+k_{4}r+k_{5}k^{2}+k_{6}r^{3}},\label{eq:rfunc_rac}
\end{equation}
where $\mathbf{k}=(k_{1},k_{2},\dots,k_{6})$ is the vector of model
parameters. This definition accommodates several models already proposed
in the literature~\cite{Ma04}. However, we will see that the shape
optimization procedure holds for any rational function.

\subsection{Extrapolation issues of radial distortion calibration}

\begin{figure}[t]
\begin{centering}
\subfloat{\centering{}%
\begin{tabular}{c}
\includegraphics[width=0.45\textwidth]{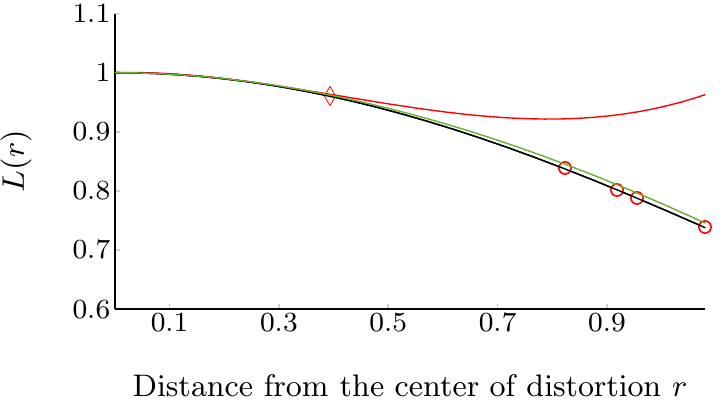}\tabularnewline
(a)\tabularnewline
\end{tabular}}\subfloat{\centering{}%
\begin{tabular}{c}
\includegraphics[width=0.45\textwidth]{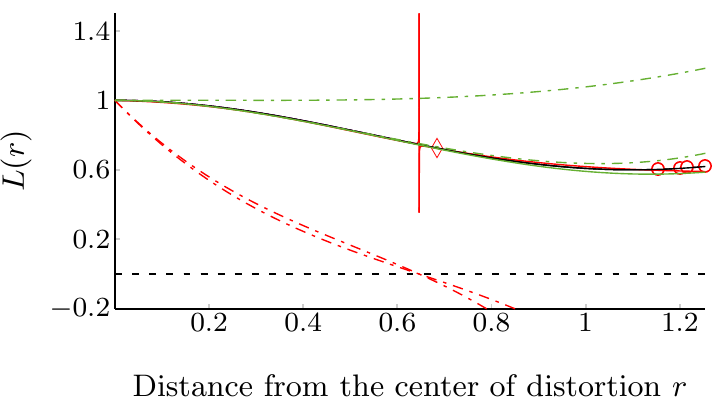}\tabularnewline
(c)\tabularnewline
\end{tabular}}\vspace{-0.3cm}

\par\end{centering}

\centering{}\subfloat{\centering{}%
\begin{tabular}{ccc}
\includegraphics[width=0.45\textwidth]{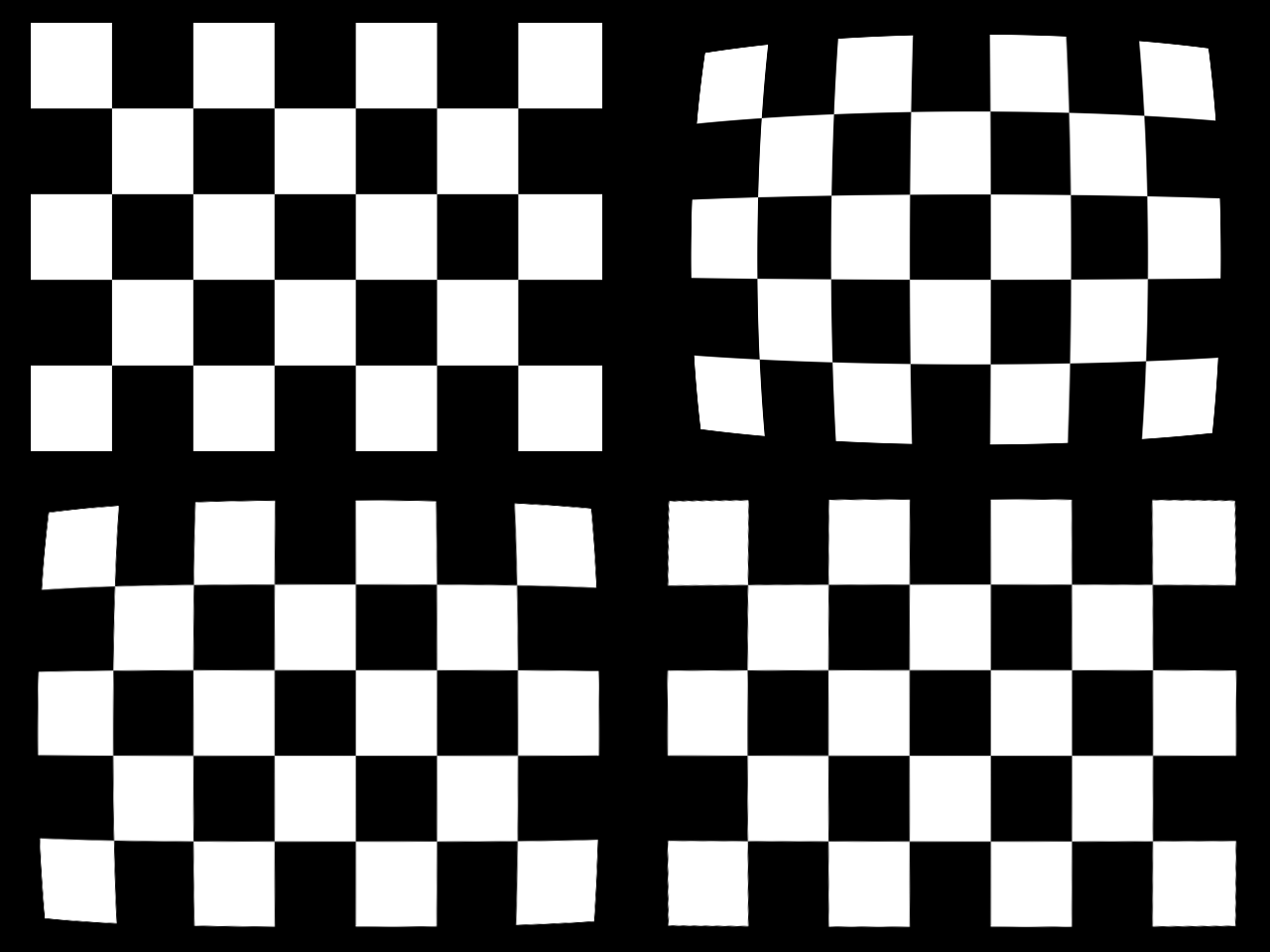} &  &
\includegraphics[width=0.45\textwidth]{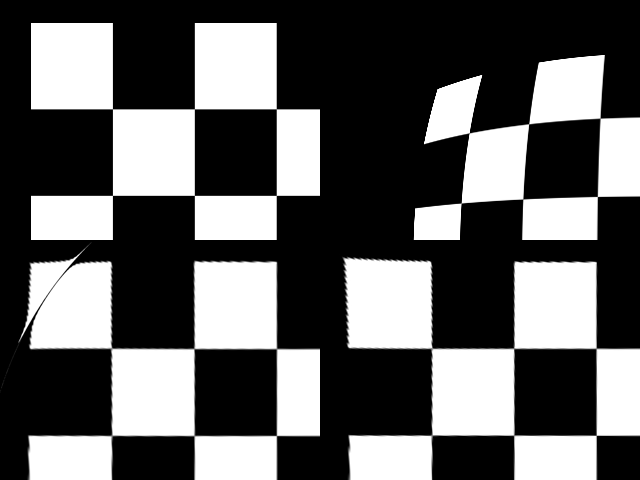}\tabularnewline
(b) &  & (d)\tabularnewline
\end{tabular}}\caption{\label{fig:distortion_issues}\emph{Calibration issues}. Examples
of issues arising while using polynomial and rational function for
radial distortion calibration. See text for details.}
\end{figure}

Let us motivate the need for the radial distortion shape optimization
by demonstrating two examples of extrapolation issues arising while
using polynomial and rational distortion models.

First, let's suppose a calibration scenario where images of a calibration
target were taken, but the image projections of the known 3D points
lie close to the center of the images with no points covering the
corners of the images. Figure~\ref{fig:distortion_issues}(a) shows
in black the graph of the amount of barrel distortion introduced by
the camera lens as a function of the distance from the center of the
radial distortion. When a polynomial distortion model $L(r)=f(r)$
is used, see Equation~\ref{eq:rfunc_rac}, in combination with an
unconstrained calibration method~\cite{Zhang00,OpenCV} (in red),
the real distortion is fitted successfully near the center of the
image on intervals where the data points are available (left of the
diamond symbol). However, the recovered polynomial quickly drifts
away elsewhere (red circles depict the distances of the projections
of the image corners). In green, a polynomial recovered by the method
proposed in this paper is shown. Here, the negativity of the first
and the second derivatives of the polynomial on the whole field of
view was enforced. This caused the model to fit the original distortion
much closer on the whole field of view. Figure~\ref{fig:distortion_issues}(b)
shows a synthetic checkerboard image (the upper left corner) and the
same image distorted by the original barrel distortion (the upper
right corner). In the lower left corner, the image is undistorted
back using the polynomial recovered by~\cite{Zhang00,OpenCV}. In
the lower right corner, the image successfully undistorted by the
polynomial recovered using the proposed shape optimization method
is shown.

Let us consider a similar calibration scenario to the one from the
previous paragraph, this time with a lens causing a mustache type
radial distortion, see Figure~\ref{fig:distortion_issues}(c). If
the radial distortion model is used, $L(r)=\frac{f(r)}{g(r)}$, the
classical calibration approach~\cite{Zhang00,OpenCV} is able to
correctly recover the original shape. However, the polynomials $f(r)$
and $g(r)$ share a common root (red dash-dot lines), which causes
a numerical instability presented as a sharp spike in $L(r)$ around
the common root---an issue we will call the \emph{zero-crossing problem}.
When the nonnegativity of $g(r)$ is enforced using the proposed approach,
not only is the correct shape recovered, but since there is now no
root in the field of view interval (green dash-dot lines), the spike
in $L(r)$ is also gone. Figure~\ref{fig:distortion_issues}(d) shows
a similar arrangement as Figure~\ref{fig:distortion_issues}(b),
now with only the upper left part of the checkerboard shown. The numerical
instability of $L(r)$ is presented as a notable ringing in the upper
left corner of the checkerboard. One can argue that the common root
is a consequence of the fact that the degrees of $f(r),g(r)$ are
higher that needed and that a model with fewer coefficients should
be used. This may be true in some cases, however, we observed just
as many situations where the lower degree polynomials resolved the
zero-crossing problem only at the cost of a considerably higher reprojection
error.

\section{Polynomials and PMI Programming\label{sec:theory}}

In this section, we present a minimal theoretical background needed
for the proposed shape optimization procedure.

\subsection{Polynomials and polynomial matrices}

An univariate polynomial $p(x)\in\Rxn n$ of degree $n\in\mathbb{N}$
is a real function defined as
\[
p(x)=p_{n}x^{n}+p_{n-1}x^{n-1}+\cdots+p_{1}x+p_{0}=\mathbf{p}^{\top}\Pxn n,
\]
where $\mathbf{p}=(p_{0},p_{1},\dots,p_{n})^{\top}\in\bRn{n+1}$ is
the vector of coefficients with a nonvanishing coefficient $p_{n}$
and $\Pxn n=(1,x,x^{2},\dots,x^{n})^{\top}$ is the canonical basis.
Let $q(x)\in\Rxn{2n}$. A~symmetric matrix $\mathtt{Q}\in\bRn{n^{\prime}\times n^{\prime}}$,
$\mathtt{Q}=(q_{i,j})$, where $n^{\prime}=n+1$, is called \emph{Gram
matrix} associated with $q(x)$ and the basis $\Pxn n$~\cite{Choi95}
if 
\begin{equation}
q(x)=\Pxnt n\,\mathtt{Q}\,\Pxn n.\label{eq:gram_matrix}
\end{equation}
Generally, there is more than one Gram matrix associated with a polynomial
$q(x)$ and we will denote the set of such matrices as $\mathcal{G}(q(x))$.The
polynomial $q(x)$ can be expressed in the elements of $\mathtt{Q}$
by simply expanding the right hand side of Equation~\ref{eq:gram_matrix}
and by comparing the coefficients. 

{}

Let $\mathbf{x}=(x_{1},x_{2},\dots,x_{d})\in\bRn d$ be a real vector
and $\ba=(\alpha_{1},\alpha_{1},\dots,\alpha_{d})\in\mathbb{N}^{d}$
an integer vector. A \emph{monomial} of degree $n=\sum\alpha_{i}$
is defined as $\mathbf{x}^{\ba}=\prod_{i=1}^{n}x_{i}^{\alpha_{i}}.$
A multivariate polynomial $p(\mathbf{x})\in\Rxxn n$ of degree $n\in\mathbb{N}$
is a mapping from $\bRn d$ to $\bR$ defined as a linear combination
of monomials up to degree $n$,
\[
p(\mathbf{x})=\sum_{|\ba|\leq n}p_{\ba}\mathbf{x}^{\ba}=\sum_{|\ba|\leq n}p_{\ba}x_{1}^{\alpha_{1}}x_{2}^{\alpha_{2}}\cdots x_{d}^{\alpha_{d}}=(p_{\ba})_{|\ba|\leq n}^{\top}(\mathbf{x}_{\ba})_{|\ba|\leq d}=\mathbf{p}^{\top}\Pxxn n,
\]
where $\mathbf{p}\in\bRn m$ is the vector of coefficients and $\Pxxn n$
is the canonical basis of $m=\tbinom{d+n}{d}$ monomials up to degree
$n$. By a polynomial matrix we will understand a symmetric matrix
whose elements are polynomials. In the next, $\mathbb{S}^{n}(\Rxxn{})$
will denote the set of $n\times n$ symmetric polynomial matrices.
The degree of $\mathtt{P}=(p_{i,j}(\mathbf{x}))\in\mathbb{S}^{n}(\Rxxn{})$
is the largest degree of all the polynomial elements of $\mathtt{P}$,
$\deg\mathtt{P}=\max_{i,j}\deg p_{i,j}(\mathbf{x})$.

Besides parameterizing polynomials by the associated Gram matrices,
we will also need to ``linearize'' them, \emph{i.e.}, to substitute
every monomial $\mathbf{x}^{\ba}$ by a new variable $y_{\ba}\in\bR$.
To do this, we define the \emph{Riesz functional} $\ell_{\mathbf{y}}\colon\mathbb{R}_{n}[\mathbf{x}]\rightarrow\mathbb{R}[\mathbf{y}]$,
a linear functional that for a $d$-variate polynomial of degree $n$,
$p(\mathbf{x})=\sum_{\ba}p_{\ba}\mathbf{x}^{\ba}$, returns an $m$-variate
polynomial of degree one, $\ell_{\mathbf{y}}(p(\mathbf{x}))=\sum_{\ba}p_{\ba}y_{\ba}$,
$m=\tbinom{d+n}{d}$. With a slight abuse of notation, we will also
use $\ell_{\mathbf{y}}$ as a matrix operator acting on $\mathbb{S}^{n}(\Rxxn{})$:
if $\mathtt{P}\in\mathbb{S}^{n}(\Rxxn{})$, then $\mathtt{P}'=\ell_{\mathbf{y}}(\mathtt{P})$
if and only if $p_{i,j}^{\prime}(\mathbf{y})=\ell_{\mathbf{y}}(p_{i,j}(\mathbf{x}))$.

\subsection{Polynomials positive on finite intervals\label{sub:positive_polynomials}}

 The shape optimization procedure presented in this paper is based
on enforcing nonnegativity of certain polynomials. Since most of the
real cameras have limited fields of view, we only need to control
the behavior of $L(r)$ for values $r\in[0,\bar{r}]$, where $\bar{r}$
is the maximal distance between the center of the radial distortion
and an (undistorted) image point. For this, we need to characterize
the set of univariate polynomials nonnegative on finite intervals.
In~\cite{Nesterov00}, based on Markov-Lukacs theorem, Nesterov showed
how to characterize such a set using positive semidefinite Gram matrices:
\begin{theorem}
\label{thm:markov-lukacs}Let $\alpha<\beta$, $p(x)\in\Rxn{}$ and
$\deg p(x)=2n$. Then $p(x)\geq0$ for all $x\in[\alpha,\beta]$ if
and only if 
\[
p(x)=s(x)+(x-\alpha)(\beta-x)t(x),
\]
where $s(x)=\Pxnt n\,\mathtt{S}\,\Pxn n$, $t(x)=\Pxnt{n-1}\,\mathtt{T}\,\Pxn{n-1}$,
such that $\mathtt{S},\mathtt{T}\succeq0$ (i.e., $\mathtt{S}\in\mathcal{G}(s(x))$,
$\mathtt{T}\in\mathcal{G}(t(x))$ are positive semidefinite Gram matrices
of polynomials $s(x)$ and $t(x)$, respectively). 

If $\deg p(x)=2n+1$, then $p(x)\geq0$ for all $x\in[\alpha,\beta]$
if and only if
\[
p(x)=(x-\alpha)s(x)+(\beta-x)t(x),
\]
where $s(x)=\Pxnt n\,\mathtt{S}\,\Pxn n$, $t(x)=\Pxnt n\,\mathtt{T}\,\Pxn n$,
such that $\mathtt{S},\mathtt{T}\succeq0$.
\end{theorem}
Even though Theorem~\ref{thm:markov-lukacs} is an equivalence, we
will only use it as an implication: as long as we will have matrices
$\mathtt{S},\mathtt{T}$ that are positive semidefinitive, Theorem~\ref{thm:markov-lukacs}
guarantees that a polynomial $p(x)$ constructed using these matrices
will be nonnegative on a given interval.

\subsection{Polynomial Matrix Inequalities\label{sub:Polynomial-Matrix-Inequalities}}

According to Theorem~\ref{thm:markov-lukacs}, a polynomial is nonnegative
on an interval as long the matrices $\mathtt{S},\mathtt{T}$ are positive
semidefinite. By combining these constraints with a polynomial cost
function, we get a problem of polynomial matrix inequalities (PMI)
programming. A PMI program can be formally defined as follows:

\begin{problem}[Polynomial matrix inequalities program]
\[
\begin{array}{rl}
\textrm{minimize} & p(\mathbf{x})\\
\textrm{subject to} & \mathtt{G}_{i}(\mathbf{x})\succeq0,\, i=1,\dots,m,\\
\textrm{where} & p(\mathbf{x})\in\mathbb{R}[\mathbf{x}],\mathtt{G}_{i}\in\mathbb{S}^{n_{i}}(\Rxxn{}).
\end{array}
\]
\label{prob:pmi}\vspace{-0.4cm}\end{problem}In general, Problem~\ref{prob:pmi}
is a hard non-convex problem. Note however, that if the cost function
$p(\mathbf{x})$ and the matrices $\mathtt{G}_{i}(\mathbf{x})$, $i=1,\dots,m$
have degree one, then Problem~\ref{prob:pmi} reduces to a linear
matrix inequality (LMI) program and as such is a semidefinite program
(SDP) solvable by any available SDP solver. In fact, most of the time
the shape optimization problems in this paper lead to such a program. 

Sometimes still, $\mathtt{G}_{i}(\mathbf{x})$ will not be linear.
In such cases, we will use the relaxation approach suggested by Henrion
and Lasserre~\cite{Henrion06}. In~\cite{Henrion06}, the authors
proposed a hierarchy of LMI programs $\mathcal{P}_{1},\mathcal{P}_{2},\dots$
that produces a monotonically non-decreasing sequence of lower bounds
$p(\mathbf{x}_{1}^{*})\leq p(\mathbf{x}_{2}^{*})\leq\dots$ on Problem~\ref{prob:pmi}
that converges to the global minimum $p(\mathbf{x}^{*})$. Practically,
the series converges to $p(\mathbf{x}^{*})$ in finitely many steps,
\emph{i.e.}, there exists $j\in\mathbb{N}$, such that $p(\mathbf{x}_{j}^{*})=p(\mathbf{x}^{*})$.
The authors also showed how this situation can be detected and how
the value of $\mathbf{x}^{*}$ can be extracted from the solution
of the relaxation by the tools of linear algebra. 

Let us show here how to construct $\mathcal{P}_{\delta}$, \emph{i.e.},
the LMI relaxation of Problem~\ref{prob:pmi} of order $\delta$;
see~\cite{Henrion06} for the technical justification of this procedure.
Let $\mathtt{G}\in\mathbb{S}^{n}(\Rxxn{})$, $n=\sum_{i=1}^{m}n_{i}$
denote a block diagonal matrix with matrices $\mathtt{G}_{i}$ on
it's diagonal. Since $(\forall i\colon\mathtt{G}_{i}(\mathbf{x})\succeq0)\Leftrightarrow\mathtt{G}(\mathbf{x})\succeq0$,
we can replace the PMI constraints $\mathtt{G}_{i}(\mathbf{x})\succeq0$
with one PMI constraint $\mathtt{G}(\mathbf{x})\succeq0$. Next, we
construct the so-called \emph{moment matrix} $\mathtt{M}_{\delta}(\mathbf{y})$
and \emph{localizing matrix} $\mathtt{M}_{\delta}(\mathtt{G},\mathbf{y})$\emph{
of }$\mathtt{G}$, defined as 
\begin{eqnarray*}
\mathtt{M}_{\delta}(\mathbf{y}) & = & \ell_{\mathbf{y}}(\Pxxn\delta\Pxxnt\delta),\\
\mathtt{M}_{\delta}(\mathtt{G},\mathbf{y}) & = & \ell_{\mathbf{y}}((\Pxxn\delta\Pxxnt\delta)\otimes\mathtt{G}),
\end{eqnarray*}
where $\otimes$ denotes the Kronecker product~\cite{Henrion06}.
Let $\gamma=1$ if $\deg\mathtt{G}\leq2$, $\gamma=\frac{\left\lceil \deg\mathtt{G}\right\rceil }{2}$
otherwise. Now, we can formally write the relaxation $\mathcal{P}_{\delta}$
as\vspace{0.0cm}\begin{problem}[LMI relaxation $\mathcal{P}_\delta$  of order $\delta$]\vspace{-0.0cm}

\[
\begin{array}{rl}
\textrm{minimize} & \ell_{\mathbf{y}}(p(\mathbf{x}))\\
\textrm{subject to} & \mathtt{M}_{\delta-\gamma}(\mathtt{G},\mathbf{y})\succeq0,\\
 & \mathtt{M}_{\delta}(\mathbf{y})\succeq0.
\end{array}
\]
\label{prob:lmi_rel}\vspace{-0.2cm}\end{problem}As the Riesz functional
$\ell_{\mathbf{y}}$ was used to ``linearize'' both the cost function
and the constraints, we can easily see that Problem~\ref{prob:lmi_rel}
is an LMI program.

\section{Shape optimization for radial distortion calibration\label{sec:shape-optimization}}

In this section, we show how to combine the results presented in Section~\ref{sec:theory}
into the radial distortion shape optimization procedure. Technically,
the procedure consists of minimization of a polynomial cost function
in the vector of radial distortion parameters $\mathbf{k}$ subject
to PMI constraints enforcing nonnegativity of certain polynomials
in the radius $r$. Such a minimization problem is a PMI program that
can be dealt with using the approach from Section~\ref{sub:Polynomial-Matrix-Inequalities}.

As mentioned in Section~\ref{sub:positive_polynomials}, we only
need to control the shape of $L(r)$ on the interval $[0,\bar{r}]$.
Note, that $\bar{r}$ is the maximal distance between the center of
the radial distortion and \emph{undistorted} image points, \emph{i.e.},
the value of $\bar{r}$ is not known prior to the actual calibration.
The value of $\bar{r}$ is therefore a user supplied parameter. Fortunately,
the proposed method is not very sensitive to the value of this parameter
and even a gross overestimate yields minima identical to the ground
truth value.

\subsection{Unconstrained radial distortion calibration}

There are several ways how to determine the vector of parameters $\mathbf{k}$
of the distortion function $L(r)$~\cite{Hartley2003b,Szeliski10}.
All we need for our shape optimization approach is a polynomial cost
function. Here, we will define and use one of such possible cost functions.
Let us rewrite Equation~\ref{eq:rdfunc} using $L(r)$ from Equation~\ref{eq:rfunc_rac}
as 

\[
g(r)\left(\begin{array}{c}
\hat{x}_{i}\\
\hat{y}_{i}
\end{array}\right)-f(r)\left(\begin{array}{c}
x_{i}\\
y_{i}
\end{array}\right)=\left(\begin{array}{c}
g(r)\,\hat{x}_{i}-f(r)\, x_{i}\\
g(r)\,\hat{y}_{i}-f(r)\, y_{i}
\end{array}\right)=\mathbf{0}.
\]
By factoring out the vector of parameters $\mathbf{k}$ and by denoting
\[
\mathtt{A}_{i}=\left(\begin{array}{cccccc}
-r\, x_{i} & -r^{2}\, x_{i} & -r^{3}\, x_{i} & \hat{x}_{i}\, r & \hat{x}_{i}\, r^{2} & \hat{x}_{i}\, r^{3}\\
-r\, y_{i} & -r^{2}\, y_{i} & -r^{3}\, y_{i} & \hat{y}_{i}\, r & \hat{y}_{i}\, r^{2} & \hat{y}_{i}\, r^{3}
\end{array}\right),\,\mathbf{b}_{i}=\left(\begin{array}{c}
x_{i}-\hat{x}_{i}\\
y_{i}-\hat{y}_{i}
\end{array}\right),
\]
we get a linear system $\mathtt{A}_{i}\mathbf{k}=\mathbf{b}_{i}$.
Now, we can stack $\mathtt{A}=(\mathtt{A}_{1}^{\top},\mathtt{A}_{2}^{\top},\dots,\mathtt{A}_{n}^{\top})^{\top}$,
$\mathbf{b}=(\mathbf{b}_{1}^{\top},\mathbf{b}_{2}^{\top},\dots\mathbf{b}_{n}^{\top})^{\top}$
and estimate the radial distortion parameters $\mathbf{k}=(k_{1},k_{2},\dots,k_{6})$
as a solution to an overdetermined system $\mathtt{A}\mathbf{k}=\mathbf{b}$
in the least square sense, \emph{i.e.}, by minimizing $\left\Vert \mathtt{A}\mathbf{k}-\mathbf{b}\right\Vert ^{2}$.
Note that for polynomial model, \emph{i.e.}, $g(x)=1$, this corresponds
to the minimization of the reprojection error. 

Let us now express the minimization of $\left\Vert \mathtt{A}\mathbf{k}-\mathbf{b}\right\Vert ^{2}$
as an LMI program. By expanding
\[
\left\Vert \mathtt{A}_{i}\mathbf{k}-\mathbf{b}_{i}\right\Vert ^{2}=(\mathtt{A}_{i}\mathbf{k}-\mathbf{b}_{i})^{\top}(\mathtt{A}_{i}\mathbf{k}-\mathbf{b}_{i})=\mathbf{k}^{\top}\mathtt{A}_{i}^{\top}\mathtt{A}_{i}\mathbf{k}-2\mathbf{b}_{i}^{\top}\mathtt{A}_{i}\mathbf{k}+\mathbf{b}_{i}^{\top}\mathbf{b}_{i}
\]
and by denoting $\mathtt{M}=\sum_{i=1}^{n}\mathtt{A}_{i}^{\top}\mathtt{A}_{i}$,
$\mathbf{m}=-2\sum_{i=1}^{n}\mathtt{A}_{i}^{\top}\mathbf{b}_{i}$,
$c=\sum_{i=1}^{n}\mathbf{b}_{i}^{\top}\mathbf{b}_{i}$, we can write
the polynomial form of the cost function as
\begin{equation}
\left\Vert \mathtt{A}\mathbf{k}-\mathbf{b}\right\Vert ^{2}=\mathbf{k}^{\top}\mathtt{M}\,\mathbf{k}+\mathbf{m}^{\top}\mathbf{k}+c.\label{eq:k_lls}
\end{equation}
As expected, Equation \ref{eq:k_lls} is a quadratic polynomial in
$\mathbf{k}$ and by construction $\mathtt{M}\succeq0$, \emph{i.e.},
$\mathtt{M}$ is a positive semidefinite matrix. Even though the cost
function is quadratic, it can be converted into a linear function
using the Schur complement trick~\cite{Boyd04}:
\[
\mathtt{F}=\left(\begin{array}{cc}
\mathtt{I} & \mathtt{L}\mathbf{k}\\
\mathbf{k}^{\top}\mathtt{L}^{\top} & -\mathbf{m}^{\top}\mathbf{k}-c+\gamma
\end{array}\right)\succeq0\,\,\Leftrightarrow\,\,\mathbf{k}^{\top}\mathtt{L}^{\top}\mathtt{L}\mathbf{k}+\mathbf{m}^{\top}\mathbf{k}+c-\gamma\leq0.
\]
By decomposing $\mathtt{M}$ as $\mathtt{M}=\mathtt{L}^{\top}\mathtt{L}$,
\emph{e.g.}, using the Cholesky or the spectral decomposition~\cite{Golub12}
(recall that $\mathtt{M}\succeq0$), we can rewrite the minimization
of Equation~\ref{eq:k_lls} as the following LMI program: \vspace{0.0cm}\begin{problem}[Unconstrained radial distortion calibration]\vspace{-0.0cm}

\[
\begin{array}{rl}
\textrm{minimize} & \gamma\\
\textrm{subject to} & \mathtt{F}=\left(\begin{array}{cc}
\mathtt{I} & \mathtt{L}\mathbf{k}\\
\mathbf{k}^{\top}\mathtt{L}^{\top} & -\mathbf{m}^{\top}\mathbf{k}-c+\gamma
\end{array}\right)\succeq0.
\end{array}
\]
\label{prob:rdistortion_calibration}\vspace{-0.2cm}\end{problem}

\subsection{Barrel distortion and the polynomial model}

As we can see from the example of barrel radial distortion in Figure~\ref{fig:distortion_issues}(a),
this type of distortion can be characterized by the negativity of
the first and the second derivatives:
\begin{equation}
\forall r\in[0,\bar{r}]\colon L'(r)\leq0\,\&\, L''(r)\leq0,\label{eq:barrel_eq_0}
\end{equation}
where $[0,\bar{r}]$ spans the field of view of the camera. If we
consider the polynomial model $L(r)=f(r)$, the constraints above
mean that we need to enforce nonnegativity of polynomials 
\[
-f'(r)=-k_{1}-2k_{2}\, r-3k_{3}\, r{}^{2},-f''(r)=-2k_{2}-6k_{3}\, r
\]
on the interval $[0,\bar{r}]$. According to Theorem~\ref{thm:markov-lukacs},
$-f'(r)\geq0$ for $\forall r\in[0,\bar{r}]$ iff
\begin{equation}
-f'(r)=-k_{1}-2k_{2}\, r-3k_{3}\, r{}^{2}=\boldsymbol{\psi}_{1}(r)^{\top}\mathtt{S}_{1}\boldsymbol{\psi}_{1}(r)+r\,(\bar{r}-r)\,\mathtt{T}_{1},\label{eq:barrel_eq_1}
\end{equation}
where 
\[
\mathtt{S}_{1}=\left(\begin{array}{cc}
s_{11} & s_{12}\\
s_{12} & s_{13}
\end{array}\right)\succeq0,\,\mathtt{T}_{1}=\left(t_{11}\right)\succeq0.
\]
By expanding the right hand side of Equation~\ref{eq:barrel_eq_1}
and by comparing the polynomial coefficients, we get a parameterization
of $\mathbf{k}$ in the elements of $\mathtt{S}_{1}$ and $\mathtt{T}_{1}$:
\begin{equation}
\left.\begin{array}{rcl}
-k_{1} & = & s_{11}\\
-2k_{2} & = & 2s_{12}+\bar{r}\, t_{11}\\
-3k_{3} & = & s_{13}-t_{11}
\end{array}\right\} \,\Rightarrow\mathbf{k}=(-s_{11},-s_{12}-{\textstyle \frac{1}{2}}\bar{r}t_{11},{\textstyle \frac{1}{3}}(t_{11}-s_{13}),0,0,0).\label{eq:barrel_k}
\end{equation}
Let's apply Theorem~\ref{thm:markov-lukacs} to $-f''(r)$ to get
the following constraint:
\begin{equation}
-f''(r)=-2k_{2}-6k_{3}\, r=r\,\mathtt{S}_{2}+(\bar{r}-r)\mathtt{T}_{2},\mathtt{S}_{2}=(s_{21})\succeq0,\mathtt{T}_{2}=(t_{21})\succeq0.\label{eq:barrel_eq_2}
\end{equation}
By combining Equations~\ref{eq:barrel_eq_2} and~\ref{eq:barrel_k},
we can express the entries of $\mathtt{S}_{2}$ and $\mathtt{T}_{2}$
in the entries of $\mathtt{S}_{1},\mathtt{T}_{1}$:
\begin{equation}
\left.\begin{array}{rcl}
-2k_{2} & = & \bar{r}\, t_{21}\phantom{\frac{1}{\bar{r}}}\\
-6k_{3} & = & s_{21}-t_{21}\phantom{\frac{1}{\bar{r}}}
\end{array}\right\} \,\Rightarrow\left\{ \begin{array}{rcl}
s_{21} & = & \frac{1}{\bar{r}}\,(2s_{12}+2\bar{r}\, s_{13}-\bar{r}t_{11})\\
t_{21} & = & \frac{2}{\bar{r}}\,(s_{12}+\frac{1}{2}\bar{r}t_{11})
\end{array}\right.
\end{equation}

Now, we have four PMI constraints on the shape of $L(r)$. If we combine
these constraints along with the parameterization of $\mathbf{k}$
from Equation~\ref{eq:barrel_k} with Problem~\ref{prob:rdistortion_calibration},
we get a radial distortion calibration problem that enforces a barrel
type distortion shape of the resulting distortion model:

\vspace{0.0cm}\begin{problem}[Barrel distortion calibration]\vspace{-0.0cm}
\[
\begin{array}{rl}
\textrm{minimize} & \gamma\\
\textrm{subject to} & \mathtt{F}\succeq0,\,\mathtt{S}_{1}\succeq0,\,\mathtt{T}_{1}=(t_{11})\succeq0,\\
 & \mathtt{S}_{2}=\left(\frac{1}{\bar{r}}\,(2s_{12}+2\bar{r}\, s_{13}-\bar{r}t_{11})\right)\succeq0,\\
 & \mathtt{T}_{2}=\left(\frac{2}{\bar{r}}\,(s_{12}+\frac{1}{2}\bar{r}t_{11})\right)\succeq0.
\end{array}
\]
\label{prob:example_1}\vspace{-0.2cm}\end{problem}Problem~\ref{prob:example_1}
is a PMI program in 5 variables $\gamma,s_{11},s_{12},s_{13},t_{11}$.
Since both the cost function and the PMI constraints have degree one,
Problem~\ref{prob:example_1} is in fact an SDP problem. Once it
is solved, the unknown distortion parameters $\mathbf{k}$ can be
easily recovered using Equation~\ref{eq:barrel_k}.

\subsection{Pincushion distortion and the division model}

Let us make an analogous analysis for the pincushion distortion shape
and the division model $L(r)=\frac{1}{g(r)}$. This type of distortion
is characterized by the nonnegativity of the first and the second
derivatives of $L(r)$ on the field of view of the camera $[0,\bar{r}]$.
From the first derivative we get the following constraint on the polynomial
denominator $g(r)$:

\[
L'(r)=\frac{-g'(r)}{g^{2}(r)}\,\,\,\,\Rightarrow\,\,\,\, L'(r)\geq0\Leftrightarrow-g'(r)\geq0.
\]
The second derivative yields a bit more complicated constraint:
\[
L''(r)=\frac{g(r)h(r)}{g^{4}(r)}=\frac{h(r)}{g^{3}(r)}\,\Rightarrow\, L''(r)\geq0\Leftrightarrow\left\{ \begin{array}{r}
(g(r)\geq0\,\&\, h(r)\geq0)\,\,\vee\,\,\,\,\,\,\,\,\,\\
\,\,\,\,\,\,\,\,\,\,\,\,(g(r)\leq0\,\&\, h(r)\leq0),
\end{array}\right.
\]
where $h(r)=2(g'(r))^{2}-g(r)g''(r)$. However, since we know that
$L(r)>0$ by definition, we only need to consider the constraints
$g(r)\geq0,h(r)\geq0$. Let us start with the constraint $g(r)\geq0$.
According to Theorem~\ref{thm:markov-lukacs}, $g(r)\geq0$ for $\forall r\in[0,\bar{r}]$
iff
\begin{equation}
g(r)=1+k_{4}r+k_{5}r^{2}+k_{6}r^{3}=\boldsymbol{\psi}_{1}(r)^{\top}\mathtt{S}_{1}\boldsymbol{\psi}_{1}(r)+(\bar{r}-r)\,\boldsymbol{\psi}_{1}(r)^{\top}\mathtt{T}_{1}\boldsymbol{\psi}_{1}(r),
\end{equation}
where
\[
\mathtt{S}_{1}=\left(\begin{array}{cc}
s_{11} & s_{12}\\
s_{12} & s_{13}
\end{array}\right)\succeq0,\,\mathtt{T}_{1}=\left(\begin{array}{cc}
t_{11} & t_{12}\\
t_{12} & t_{13}
\end{array}\right)\succeq0.
\]
This leads to the following parameterization of $\mathbf{k}$ as well
as to a constraint on the variable $t_{11}$:
\begin{equation}
\left.\begin{array}{rcl}
1 & = & \bar{r}\, t_{11}\\
k_{4} & = & s_{11}-t_{11}+2\bar{r}t_{12}\\
k_{5} & = & 2s_{12}-2t_{12}+\bar{r}t_{13}\\
k_{6} & = & s_{13}-t_{13}
\end{array}\right\} \,\Rightarrow\left\{ \begin{array}{rcl}
\mathbf{k} & = & (0,0,0,s_{11}-t_{11}+2\bar{r}t_{12},\\
 &  & \,\,\,\,\,2s_{12}-2t_{12}+\bar{r}t_{13},s_{13}-t_{13})\\
t_{11} & = & \frac{1}{\bar{r}}
\end{array}\right.\label{eq:pincushion_2}
\end{equation}
By applying Theorem~\ref{thm:markov-lukacs} to the constraint $-g'(r)\geq0$,
we get 
\begin{equation}
-g'(r)=-k_{4}-2k_{5}\, r-3k_{6}\, r{}^{2}=\boldsymbol{\psi}_{1}(r)^{\top}\mathtt{S}_{2}\boldsymbol{\psi}_{1}(r)+r\,(\bar{r}-r)\,\mathtt{T}_{2},
\end{equation}
where
\[
\mathtt{S}_{2}=\left(\begin{array}{cc}
s_{21} & s_{22}\\
s_{22} & s_{23}
\end{array}\right)\succeq0,\mathtt{T}_{2}=\left(t_{21}\right)\succeq0.
\]
As in the case of the barrel distortion optimization, we can express
the entries of $\mathtt{S}_{2}$ and $\mathtt{T}_{2}$ in the entries
of $\mathtt{S}_{1},\mathtt{T}_{1}$. This time, however, we have more
variables than equations and we have to set one of the entries free---we
chose $s_{22}$:
\begin{equation}
\left.\begin{array}{rcl}
-3k_{4} & = & s_{21}\\
-2k_{5} & = & 2s_{22}+\bar{r}t_{21}\\
-3k_{6} & = & s_{23}-t_{21}
\end{array}\right\} \,\Rightarrow\left\{ \begin{array}{rcl}
s_{21} & = & t_{11}-s_{11}-2\bar{r}t_{12}\\
s_{23} & = & -\frac{1}{\bar{r}}(s_{12}+2s_{22}-4t_{12}+\bar{r}(3s_{13}-t_{13}))\\
t_{21} & = & -\frac{1}{\bar{r}}(2s_{12}+s_{22}-2t_{12}+\bar{r}t_{13})
\end{array}\right.
\end{equation}
The final constraint is the most complicated because of the quadratic
monomials in $\mathbf{k}$: $h(r)>0$ for $\forall r\in[0,\bar{r}]$
iff
\begin{eqnarray}
h(r) & = & (6k_{6}r^{2}+4k_{5}r+2k_{4})(3k_{6}r^{2}+2k_{5}+k_{4})-\nonumber \\
 &  & \,\,\,\,\,\,\,\,\,\,\,-(2k_{5}+6k_{6}r)(k_{6}r^{3}-k_{5}r^{2}+k_{4}r+1)\label{eq:pincushion_3}\\
 & = & \boldsymbol{\psi}_{2}(r)^{\top}\mathtt{S}_{3}\boldsymbol{\psi}_{2}(r)+(\bar{r}-r)\,\boldsymbol{\psi}_{1}(r)^{\top}\mathtt{T}_{3}\boldsymbol{\psi}_{1}(r),\nonumber 
\end{eqnarray}
where
\[
\mathtt{S}_{3}=\left(\begin{array}{ccc}
s_{31} & s_{32} & s_{33}\\
s_{32} & s_{34} & s_{35}\\
s_{33} & s_{35} & s_{36}
\end{array}\right)\succeq0,\mathtt{T}_{3}=\left(\begin{array}{cc}
t_{31} & t_{32}\\
t_{32} & t_{33}
\end{array}\right)\succeq0.
\]
Equation~\ref{eq:pincushion_3} gives us 5 constraints on 9 entries
of $\mathtt{S}_{3}$ and $\mathtt{T}_{3}$. We chose to set free variables
$s_{32},s_{34},s_{36},t_{32}$; System~\ref{eq:pincushion_4} shows
the form of the remaining 5 variables. Finally, we can combine these
6 PMI constraints, Problem~\ref{prob:rdistortion_calibration} and
the parameterization of $\mathbf{k}$ from Equation~\ref{eq:pincushion_2}
into a radial distortion calibration problem that enforces a pincushion
type distortion shape:

\begin{figure*}
\addtolength{\arraycolsep}{-0.9mm}
\begin{equation}
\left.\begin{array}{rcl}
{\scriptstyle 12k_{6}^{2}} & {\scriptstyle =} & {\scriptstyle s_{36}-t_{33}}\\
{\scriptstyle 16k_{5}k_{6}} & {\scriptstyle =} & {\scriptstyle 2s_{35}-2t_{32}+\bar{r}t_{33}}\\
{\scriptstyle 6k_{5}^{2}+6k_{4}k_{6}} & {\scriptstyle =} & {\scriptstyle 2s_{33}+s_{34}-t_{31}+2\bar{r}t_{32}}\\
{\scriptstyle 6k_{4}k_{5}-6k_{6}} & {\scriptstyle =} & {\scriptstyle 2s_{32}+\bar{r}t_{31}}\\
{\scriptstyle 2k_{4}^{2}-2k_{5}} & {\scriptstyle =} & {\scriptstyle s_{31}}
\end{array}\right\} \negmedspace\Rightarrow\negmedspace\left\{ \begin{array}{rcl}
{\scriptstyle s_{31}} & {\scriptstyle =} & {\scriptstyle 4t_{12}-s_{12}-2\bar{r}t_{31}+2(s_{11}-t_{11}+2\bar{r}t_{12})^{2}}\\
{\scriptstyle s_{33}} & {\scriptstyle =} & {\scriptstyle -\frac{1}{2\bar{r}}(6s_{13}+2s_{32}-6t_{13}+\bar{r}s_{34}}\\
 &  & {\scriptstyle -6\bar{r}(2s_{12}-2t_{12}+\bar{r}t_{13})^{2}+2\bar{r}^{2}t_{32}-}\\
 &  & {\scriptstyle -6(s_{11}-t_{11}+2\bar{r}t_{12})}\\
 &  & {\scriptstyle (2s_{12}-2t_{12}+\bar{r}t_{13}+\bar{r}s_{13}-\bar{r}t_{13}))}\\
{\scriptstyle s_{35}} & {\scriptstyle =} & {\scriptstyle t_{32}-\frac{2}{\bar{r}}s_{36}{\scriptstyle +6\bar{r}(s_{13}-t_{13})^{2}}+}\\
 &  & {\scriptstyle +8(s_{13}-t_{13})(2s_{12}-2t_{12}+\bar{r}t_{13})}\\
{\scriptstyle t_{31}} & {\scriptstyle =} & {\scriptstyle -\frac{2}{\bar{r}}(3s_{13}+s_{32}-3t_{13}-}\\
 &  & {\scriptstyle 3(2s_{12}-2t_{12}+\bar{r}t_{13})(s_{11}-t_{11}+2\bar{r}t_{12}))}\\
{\scriptstyle t_{33}} & {\scriptstyle =} & {\scriptstyle s_{36}-12(s_{13}-t_{13})^{2}}
\end{array}\right.\label{eq:pincushion_4}
\end{equation}
\addtolength{\arraycolsep}{0.9mm}\hrulefill
\end{figure*}

\vspace{0.0cm}\begin{problem}[Pincushion distorion calibration]\vspace{-0.0cm}
\[
\begin{array}{rl}
\textrm{minimize} & \gamma\\
\textrm{subject to} & \mathtt{F}\succeq0,\,\mathtt{S}_{1}\succeq0,\,\mathtt{T}_{1}\succeq0,\mathtt{S}_{2}\succeq0,\,\mathtt{T}_{2}\succeq0,\mathtt{S}_{3}\succeq0,\,\mathtt{T}_{3}\succeq0.
\end{array}
\]
\label{prob:example_2}\vspace{-0.2cm}\end{problem}Problem~\ref{prob:example_2}
is a PMI program in 11 variables $\gamma$, $s_{11}$, $s_{12}$,
$s_{13}$, $t_{12}$, $t_{13}$, $s_{22}$, $s_{32},$ $s_{34}$,
$s_{36}$, and $t_{32}$. Since $\mathtt{S}_{3}$ and $\mathtt{T}_{3}$
are polynomial matrices of degree 2, Problem~\ref{prob:example_2}
has to be dealt with using the relaxation scheme from Section~\ref{sub:Polynomial-Matrix-Inequalities}.

\subsection{Zero-crossing problem of the rational model}

Also the zero-crossing problem of the rational model $L(r)=\frac{f(r)}{g(r)}$
can be dealt with using the proposed shape optimization technique.
A sufficient condition for avoiding a common root of the polynomials
$f(r)$ and $g(r)$ on the interval $[0,\bar{r}]$ is to force at
least one on them to have no root. Here, we decided on enforcing the
constraint
\begin{equation}
\forall r\in\left\langle 0,\bar{r}\right\rangle \colon g(r)-p\geq0,\textrm{ where }p>0.\label{eq:zero_crossing1}
\end{equation}
Since Theorem~\ref{thm:markov-lukacs} guarantees only nonnegativity
of a polynomial, we need a strictly positive parameter $p$ to enforce
strict positivity of $g(r)$. Even though parameter $p$ must be user
supplied, the method is not overly sensitive to its value; in our
experiments, we set $p=0.1$. By applying Theorem~\ref{thm:markov-lukacs}
to the above constraint and the interval $[0,\bar{r}]$, we get
\[
g(r)-p=1-p+k_{4}r+k_{5}r^{2}+k_{6}r^{3}=\boldsymbol{\psi}_{1}(r)^{\top}\mathtt{S}_{1}\boldsymbol{\psi}_{1}(r)+(\bar{r}-r)\,\boldsymbol{\psi}_{1}(r)^{\top}\mathtt{T}_{1}\boldsymbol{\psi}_{1}(r),
\]
where
\[
\mathtt{S}_{1}=\left(\begin{array}{cc}
s_{11} & s_{12}\\
s_{12} & s_{13}
\end{array}\right)\succeq0,\,\mathtt{T}_{1}=\left(\begin{array}{cc}
t_{11} & t_{12}\\
t_{12} & t_{13}
\end{array}\right)\succeq0.
\]
This yields a parameterization of $\mathbf{k}$ as well as a constraint
on $t_{11}$:
\begin{equation}
\left.\begin{array}{rcl}
1-p & \!\!=\!\! & \bar{r}\, t_{11}\\
k_{4} & \!\!=\!\! & s_{11}-t_{11}+2\bar{r}t_{12}\\
k_{5} & \!\!=\!\! & 2s_{12}-2t_{12}+\bar{r}t_{13}\\
k_{6} & \!\!=\!\! & s_{13}-t_{13}
\end{array}\right\} \Rightarrow\left\{ \begin{array}{rcl}
\mathbf{k} & \!\!=\!\! & (k_{1},k_{2},k_{3},s_{11}-t_{11}+2\bar{r}t_{12},\\
 &  & \,\,\,\,\,2s_{12}-2t_{12}+\bar{r}t_{13},s_{13}-t_{13})\\
t_{11} & \!\!=\!\! & \frac{1-p}{\bar{r}}
\end{array}\right.\label{eq:zero_crossing2}
\end{equation}
Again, by combining the two PMI constraints with Problem~\ref{prob:rdistortion_calibration}
and the parameterization of $\mathbf{k}$ from Equation~\ref{eq:zero_crossing2},
we get a radial distortion calibration problem that eliminates the
zero-crossing problem: \vspace{0.0cm}\begin{problem}[Zero-crossing distortion calibration]~\vspace{-0.3cm}
\[
\begin{array}{rl}
\textrm{minimize} & \gamma\\
\textrm{subject to} & \mathtt{F}\succeq0,\,\mathtt{S}_{1}\succeq0,\,\mathtt{T}_{1}\succeq0.
\end{array}
\]
\label{prob:example_3}\vspace{-0.6cm}\end{problem}Problem~\ref{prob:example_3}
is an LMI program in 9 variables $\gamma,s_{11},s_{12},s_{13},t_{12},t_{13},k_{1},k_{2},k_{3}$.

\subsection{Shape optimization in Camera Calibration Procedure\label{sec:calibration}}

All of the calibration problems presented in this paper expect the
projection coordinates $x_{i}$, $y_{i}$, $\hat{x}_{i}$, and $\hat{y}_{i}$
to be known, see Equation~\ref{eq:rdfunc}. This assumes a known
calibration target $\mathbf{X}_{i}\in\bR^{3}$ as well as known camera
parameters $\mathtt{R}\in SO(3)$, $\mathbf{t}\in\bR^{3}$, and the
calibration matrix $\mathtt{K}\in\bR^{3\times3}$.  A straightforward
idea how to fold the shape optimized radial distortion calibration
into the camera calibration procedure is to first perform ``classical''
camera calibration~\cite{Tsai86,Zhang00,Hartley07-2,Tardif09}, including
radial distortion estimation. Once the projection coordinates are
known, the shape optimized radial distortion calibration can be performed
to replace the radial distortion parameters estimated by a classical
method. One might argue that the quality of such a solution could
be compromised, since different error functions may be considered
by the camera and the shape optimization calibration methods. To mitigate
this problem, we suggest an alternating approach to ``shape-optimize''
the results of the classical camera calibration: first, the shape
optimization procedure is performed, followed by a bundle adjustment~\cite{Triggs00}
step where the radial distortion parameters are fixed. This can be
repeated in a loop for a fixed number of times, or until desired convergence
is reached.

\section{Experiments\label{sec:experiments}}

\begin{center}
\vspace{-0.8cm}

\par\end{center}

To validate the proposed approach, this section presents several experimental
results on synthetic as well as real world datasets. We implemented
Problems~\ref{prob:example_1}, \ref{prob:example_2}, and \ref{prob:example_3}
in \textsc{Matlab} using Yalmip toolbox~\cite{yalmip04} with SeDuMi~\cite{Sturm99}
as the underlying SDP solver. Yalmip toolbox is a modeling language
that can be used to solve LMI as well as PMI programs, which it automatically
translates into LMI relaxations using the scheme presented in Section~\ref{sub:Polynomial-Matrix-Inequalities}.
All of the resulting SDP programs were solved under a second on an
Intel i7 3.50GHz based desktop computer running Linux and 64bit \textsc{Matlab.}

\begin{figure}[t]
\centering{}%
\begin{tabular}{ccccc}
\includegraphics[width=0.3\textwidth]{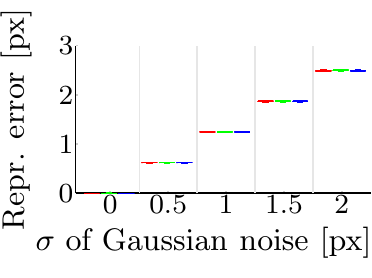} &  &
\includegraphics[width=0.3\textwidth]{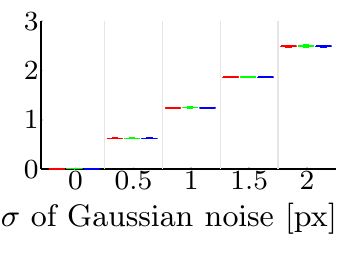} &  &
\includegraphics[width=0.3\textwidth]{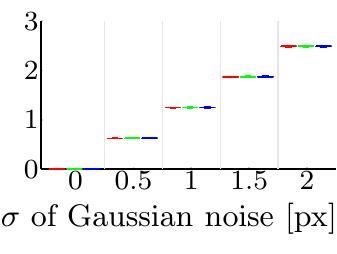}\tabularnewline
 &  & \vspace{-0.6cm}
 &  & \tabularnewline
(a) & \hspace{0.1cm} & (c) & \hspace{0.1cm} & (e)\tabularnewline
\includegraphics[width=0.3\textwidth]{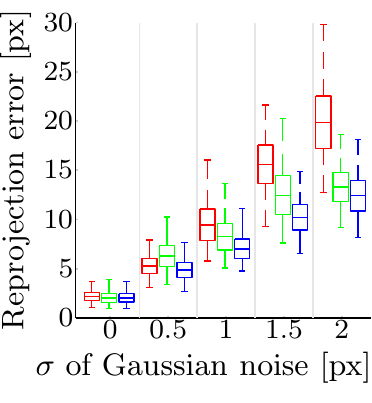} &  &
\includegraphics[width=0.3\textwidth]{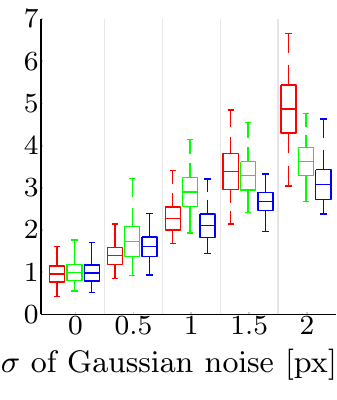} &  &
\includegraphics[width=0.3\textwidth]{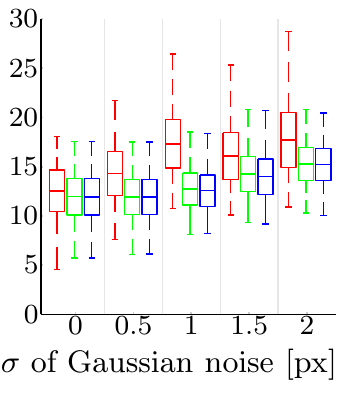}\tabularnewline
 &  & \vspace{-0.6cm}
 &  & \tabularnewline
(b) &  & (d) &  & (f)\tabularnewline
\end{tabular}\caption{\emph{Image noise experiment}. Methods BA, SO, and ASO in red, green,
and blue, respectively, on calibration and validation data point sets.
(a--b) barrel distortion, (c--d) pincushion distortion (e--f) zero-crossing
problem.\label{fig:synth_exp}}
\vspace{-0.8cm}
\end{figure}

\textbf{Synthetic experiment}. In the synthetic experiment, we studied
the performance of the proposed method with respect to the image noise.
We generated a synthetic $16\!\times\!16$ planar calibration target.
A scene consisted of 9 random $640\!\times\!480$~pixel cameras randomly
positioned on a hemisphere around the target and rotated to face its
center. The focal length was set to approx. $540$ px and the distances
of the camera centers from the target were set up so that the target
(calibration data point set) covered only the middle part of the field
of view, approx $50\%$. For each of the three model-shape problem
combinations, we generated 100 scenes and corrupted the projections
of the calibration target by an increasing amount of Gaussian image
noise in 5 levels, standard deviation $\sigma\in[0,2]$ px in $\nicefrac{1}{2}\,\textrm{px}$
steps. We calibrated all scenes with OpenCV~\cite{OpenCV} made to
disregard the radial distortion component. We compare three methods:
the first method (BA) is the bundle adjustment method that included
the respective radial distortion model performed together with the
OpenCV calibration results, the second method (SO) is the respective
shape-optimization method performed after the BA step, and the last
method (ASO) is the alternating approach from Section~\ref{sec:calibration},
fixed to 10 iterations. 

\emph{Barrel distortion}. First, we experimented with the barrel distortion
and the polynomial model $L(r)=f(r)$. Figure~\ref{fig:synth_exp}(a)
shows the mean of the reprojection errors on the calibration data
point set for methods BA, SO, and ASO using \textsc{Matlab} function
\texttt{boxplot}. The methods show identical performance, however
when a validation data set of points covering the whole field of view
is used, see Figure~\ref{fig:synth_exp}(a), we see both SO and ASO
outperforming the classical BA approach.

\emph{Pincushion distortion}. Next, Figures~\ref{fig:synth_exp}(c--d)
show the analogous measure for the pincushion distortion and the division
model $L(r)=\frac{1}{g(x)}$. Here, both BA and shape-optimization
methods perform significantly better on the validation data point
set. Still, we can see superior performance of SO and ASO as the noise
increases.

\emph{Zero-crossing}. Finally, we experimented with the rational model
$L(r)=\frac{f(r)}{g(r)}$ and the mustache type distortion. Figure~\ref{fig:synth_exp}(d)
shows identical performance on the calibration dataset. On the other
hand, we can see poor performance on the validation data point set
even if no noise is present, Figure~\ref{fig:synth_exp}(e). This
is caused by the fact that too few calibration points were on the
outer parts of the field of view where the convexity of the distortion
function changes. Again, we see better performance of SO and ASO methods.

\begin{figure}[t]
\centering{}\subfloat[]{\begin{centering}
\includegraphics[width=0.4\textwidth]{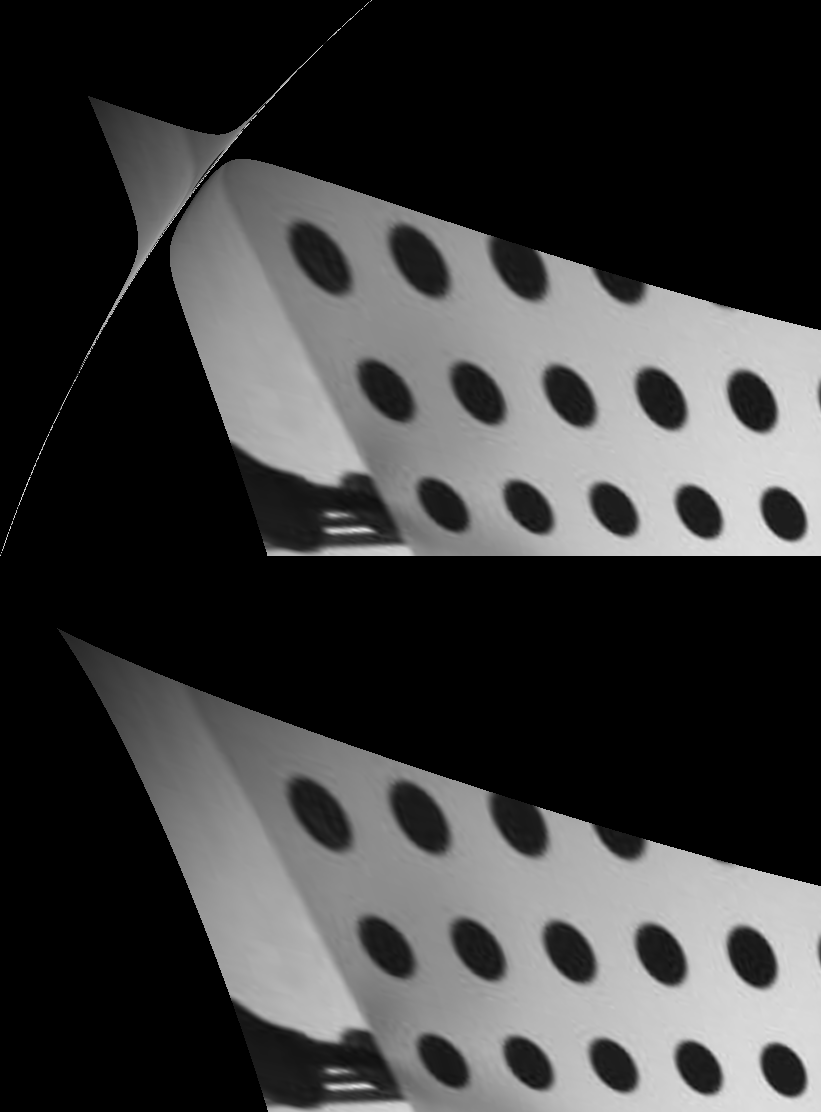}\vspace{0.5cm}

\par\end{centering}

\begin{centering}

\par\end{centering}

}\hspace{0.6cm}\subfloat[]{\centering{}\includegraphics[width=0.5\textwidth]{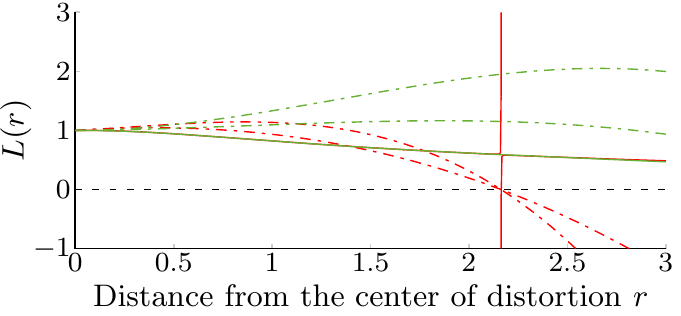}}\caption{\emph{Real experiment}. Correction of the zero-crossing problem of
the rational model.\label{fig:real_exp}}
\vspace{-0.6cm}
\end{figure}

\textbf{Real experiment}. In the real experiment, we calibrated a
2~MPix camera from Point Grey's Ladybug 3 system~\cite{Ladybug}
using 12 images of a known $28\!\times\!20$ planar target. Calibration
using BA method and the rational model introduced quite noticeable
zero-crossing problem. As expected, calibration using ASO method does
not suffer from this type of problem. In this experiment, we set $\bar{r}=4$
and $p=0.1$. Figure~\ref{fig:real_exp}(a) shows the upper left
corner of a rectified calibration image using $\mathbf{k}$ provided
by methods BA and ASO, respectively. Figure~\ref{fig:real_exp}(b)
shows the shape of the BA calibration function in red and the ASO
calibration in green.

\section{Conclusion}

\begin{center}
\vspace{-0.8cm}

\par\end{center}

The aim of this work was not to argue for a specific radial distortion
model, but to point out extrapolation problems inherent to all polynomial
and rational models. We solved these problems by enforcing a predetermined
shape of the distortion function. For most shapes and models, the
proposed approach leads to small semidefinite programming problems
that can be solved fast and globally optimally. We also showed how
to deal with shapes and models that lead to PMI problems using a LMI
relaxation scheme. We showed experimentally that in terms of the reprojection
error on the known data points the proposed approach provides radial
distortion models that are equivalent to those provided by the classical
bundle adjustment approach, yet with the added value of having the
correct shape that mollifies or completely removes all extrapolation
issues.


\begin{thebibliography}{10}

\bibitem{Ladybug}
Ladybug 3 camera.
\newblock{\tt www.ptgrey.com/products/ladybug3}.

\bibitem{OpenCV}
Open source computer vision library.
\newblock {\tt www.opencv.org}.

\bibitem{Boyd04}
Stephen Boyd and Lieven Vandenberghe.
\newblock {\em Convex Optimization}.
\newblock Cambridge University Press, March 2004.

\bibitem{Brown66}
Duane~C. Brown.
\newblock Decentering distortion of lenses.
\newblock {\em Photometric Engineering}, 32(3):444--462, 1966.

\bibitem{Brown71}
Duane~C. Brown.
\newblock Close-range camera calibration.
\newblock {\em Photogrammetric Engineering}, 37(8):855--866, 1971.

\bibitem{Choi95}
Man-Duen Choi, Tsit~Yuen Lam and Bruce Reznick.
\newblock Sums of squares of real polynomials.
\newblock In {\em Proceedings of Symposia in Pure mathematics}, volume~58,
  pages 103--126. American Mathematical Society, 1995.

\bibitem{Golub12}
Gene~H. Golub and Charles~F. Van~Loan.
\newblock {\em Matrix computations}, volume~3.
\newblock Johns Hopkins University Press, 2012.

\bibitem{Hartley07-2}
Richard Hartley and Sing~Bing Kang.
\newblock Parameter-free radial distortion correction with center of distortion
  estimation.
\newblock {\em Pattern Analysis and Machine Intelligence, IEEE Transactions
  on}, 29(8):1309--1321, 2007.

\bibitem{Hartley2003b}
Richard Hartley and Andrew Zisserman.
\newblock {\em Multiple view geometry in computer vision}.
\newblock Cambridge University, Cambridge, 2nd edition, 2003.

\bibitem{Henrion06}
Didier Henrion and Jean-Bernard Lasserre.
\newblock Convergent relaxations of polynomial matrix inequalities and static
  output feedback.
\newblock {\em IEEE Transactions on Automatic Control}, 51(2):192--202, 2006.

\bibitem{yalmip04}
Johan~L\"{o}fberg.
\newblock {YALMIP}: A toolbox for modeling and optimization in {MATLAB}.
\newblock In {\em Proceedings of the IEEE Symposium on CACSD}, Taipei, Taiwan, 2004.

\bibitem{Ma04}
Lili Ma, YangQuan Chen, and Kevin~L. Moore.
\newblock Rational radial distortion models of camera lenses with analytical
  solution for distortion correction.
\newblock {\em International Journal of Information Acquisition},
  1(02):135--147, 2004.

\bibitem{Nesterov00}
Yurii Nesterov.
\newblock Squared functional systems and optimization problems.
\newblock In {\em High performance optimization}, pages 405--440. Springer,
  2000.

\bibitem{Slama80}
Chester~C. Slama, Charles Theurer, Soren~W. Henriksen, et~al.
\newblock {\em Manual of photogrammetry.}
\newblock Number Ed. 4. American Society of photogrammetry, 1980.

\bibitem{Sturm99}
Jos F. Sturm.
\newblock Using {SeDuMi} 1.02, a {MATLAB} toolbox for optimization over
  symmetric cones.
\newblock {\em Optimization Methods and Software}, 11--12:625--653, 1999.

\bibitem{Sturm11}
Peter Sturm, Srikumar Ramalingam, Jean-Philippe Tardif, Simone Gasparini, and
  Joao Barreto.
\newblock Camera models and fundamental concepts used in geometric computer
  vision.
\newblock {\em Foundations and Trends in Computer Graphics and Vision},
  6(1--2):1--183, 2011.

\bibitem{Szeliski10}
Richard Szeliski.
\newblock {\em Computer vision: algorithms and applications}.
\newblock Springer, 2010.

\bibitem{Tardif09}
Jean-Philippe Tardif, Peter Sturm, Martin Trudeau, and Sebastien Roy.
\newblock Calibration of cameras with radially symmetric distortion.
\newblock {\em IEEE Transactions
  on Pattern Analysis and Machine Intelligence}, 31(9):1552--1566, 2009.

\bibitem{Triggs00}
Bill Triggs, Philip~F. McLauchlan, Richard~I. Hartley, and Andrew~W.
  Fitzgibbon.
\newblock Bundle adjustment - a modern synthesis.
\newblock In {\em ICCV '99: Proceedings of the International Workshop on Vision
  Algorithms}, pages 298--372, London, UK, 2000.

\bibitem{Tsai86}
Roger~Y. Tsai.
\newblock An efficient and accurate camera calibration technique for 3d machine
  vision.
\newblock In {\em Proc. IEEE Conf. on Computer Vision and Pattern Recognition}, 1986.

\bibitem{Zhang00}
Zhengyou Zhang.
\newblock A flexible new technique for camera calibration.
\newblock {\em IEEE Transactions
  on Pattern Analysis and Machine Intelligence}, 22(11):1330--1334, 2000.

\end{thebibliography}
\end{document}